\documentclass[11pt]{article}
\usepackage{epsf, graphicx}
\newfont{\Bb}{msbm10 scaled\magstep1}
\newfont{\Bbs}{msbm10 scaled\magstep0}

\newcommand{\fq}{\mbox{\Bb F}_{q}}
\newcommand{\fqk}{\mbox{\Bb F}_{q^{k}}}

\newcommand{\fp}{\mbox{\Bb F}_{p}}
\newcommand{\fqbar}{\bar{\mbox{\Bb F}}_{q}}

\newcommand{\C}{\mbox{\Bb C}}
\newcommand{\G}{\mbox{\Bb G}}
\newcommand{\PP}{\mbox{\Bb P}}
\newcommand{\A}{\mbox{\Bb A}}

\newcommand{\Q}{\mbox{\Bb Q}}
\newcommand{\Z}{\mbox{\Bb Z}}

\newcommand{\R}{\mbox{\Bb R}}
\newtheorem{theorem}{Theorem}[section]
\newtheorem{corollary}[theorem]{Corollary}
\newtheorem{lemma}[theorem]{Lemma}

\newtheorem{proposition}[theorem]{Proposition}
\newtheorem{definition}[theorem]{Definition}
\newtheorem{conjecture}[theorem]{Conjecture}
\newtheorem{question}[theorem]{Question}

\newenvironment{exafont}{\begin{bf}}{\end{bf}}

\begin{document}

\title{Mirror Symmetry For Zeta Functions}

\author{Daqing Wan\footnote{Partially supported by NSF. It is a pleasure to thank 
V. Batyrev, P. Candelas, H. Esnault, K.F. Liu, Y. Ruan, S.T. Yau for   
helpful discussions. The paper is motivated 
by some open questions in my lectures at the 2004 Arizona Winter School.}}

\maketitle
\centerline{Institute of Mathematics, Chinese Academy of Sciences, Beijing, P.R. China}
\centerline{Department of Mathematics, University of California, Irvine, CA 92697-3875}
\centerline{dwan@math.uci.edu}

\begin{abstract}
In this paper, we study the relation 
between the zeta function of a Calabi-Yau hypersurface and the zeta 
function of its mirror. Two types of arithmetic relations 
are discovered. This motivates us to formulate two general arithmetic 
mirror conjectures for the zeta functions of a mirror pair of 
Calabi-Yau manifolds.  

\end{abstract}

\tableofcontents

\section{Introduction}

In this section, we describe two mirror relations between 
the zeta function of a Calabi-Yau hypersurface in a projective space 
and the zeta function of its mirror manifold. Along the way, we 
make comments and conjectures about what to expect in the 
general case.  

Let $d$ be a positive integer. 
Let $X$ and $Y$ be two $d$-dimensional smooth projective Calabi-Yau varieties  
over $\C$. A necessary condition (the topological mirror test) for $X$ and $Y$ 
to be a mirror pair is that their Hodge numbers satisfy the Hodge symmetry: 
\begin{equation}
h^{i,j}(X) = h^{d-i, j}(Y), \ 0\leq i, j \leq d. 
\end{equation}
In particular, their Euler characteristics are related by 
\begin{equation}
e(X) = (-1)^d e(Y).
\end{equation}
In general, there is no known rigorous algebraic geometric definition 
for a mirror pair, although many examples of mirror pairs are known 
at least conjecturally. Furthermore, it does not make sense to speak 
of ``the mirror'' of $X$ as the mirror variety usually comes in a family.  
In some cases, the mirror does not exist. This is the case for 
rigid Calabi-Yau $3$-fold $X$, since the rigid condition $h^{2,1}(X)=0$ would 
imply that $h^{1,1}(Y)=0$ which is impossible. 

We shall assume  that $X$ and $Y$ are a given 
mirror pair in some sense and are defined over a number field or a finite 
field. We are interested in how the zeta function of $X$ is related to the 
zeta function of $Y$. 
Since there is no algebraic geometric definition for 
$X$ and $Y$ to be a mirror pair, it is difficult to study the possible 
symmetry between their zeta functions in full generality. 
On the other hand, there do have many explicit examples 
and constructions which at least conjecturally give a mirror pair, 
most notably in the toric hypersurface setting as constructed by Batyrev \cite{Ba1}. 
Thus, we shall first examine an explicit example 
and see what kind of relations can be proved for their zeta functions in this case. 
This would then suggest what to expect in general. 

Let $n\geq 2$ be a positive integer.
We consider the universal family of Calabi-Yau complex hypersurfaces of degree $n+1$ 
in the projective space $\PP^n$. Its mirror family is a one parameter family of toric 
hypersurfaces. To construct the mirror family, we consider 
the one parameter subfamily $X_{\lambda}$ of complex projective hypersurfaces 
of degree $n+1$ in $\PP^n$ 
defined by 
$$f(x_1,\cdots, x_{n+1})= x_1^{n+1}+\cdots +x_{n+1}^{n+1}+\lambda x_1\cdots x_{n+1}=0,$$
where $\lambda\in \C$ is the parameter. 
The variety $X_{\lambda}$ is a Calabi-Yau manifold 
when $X_{\lambda}$ is smooth. 
Let $\mu_{n+1}$ denote the group of 
$(n+1)$-th roots of unity. Let 
$$G=\{(\zeta_1,\cdots, \zeta_{n+1})|\zeta_i^{n+1}=1, \zeta_1\cdots \zeta_{n+1}=1\}/\mu_{n+1} 
\cong (\Z/(n+1)\Z)^{n-1},$$
where $\mu_{n+1}$ is embedded in $G$ via the diagonal embedding. 
The finite group $G$ acts on $X_{\lambda}$ by 
$$(\zeta_1,\cdots, \zeta_{n+1})(x_1,\cdots, x_{n+1})=(\zeta_1x_1, \cdots, \zeta_{n+1}x_{n+1}).$$
The quotient $X_{\lambda}/G$ is a projective toric hypersurface $Y_{\lambda}$ in the 
toric variety $\PP_{\Delta}$, where $\PP_{\Delta}$ is the simplex in $\R^n$ with 
vertices $\{ e_1,\cdots, e_n, -(e_1+\cdots e_n)\}$ and the $e_i$'s are the standard 
coordinate vectors in $\R^n$. Explicitly, the variety $Y_{\lambda}$ is 
the projective closure in $\PP_{\Delta}$ of the affine toric hypersurface in $\G_m^n$  
defined by 
$$g(x_1,\cdots, x_{n})= x_1+\cdots +x_n+{1\over x_1\cdots x_{n}}+\lambda=0.$$

Assume that $X_{\lambda}$ is smooth. Then, $Y_{\lambda}$ is a (singular) mirror 
of $X_{\lambda}$. It is an orbifold. 
If $W_{\lambda}$ is a smooth crepant resolution 
of $Y_{\lambda}$, 
then the pair $(X_{\lambda}, W_{\lambda})$ 
is called a mirror 
pair of Calabi-Yau manifolds. 
Such a resolution  exists for this example 
but not unique if $n\geq 3$. The number of rational points and the zeta function 
are independent of the choice of the crepant resolution.  
We are interested in understanding how the arithmetic 
of $X_{\lambda}$ is related to the arithmetic of $W_{\lambda}$, in particular how 
the zeta function of $X_{\lambda}$ is related to the zeta function of $W_{\lambda}$. 
Our main concern in this paper is to consider Calabi-Yau manifolds 
over finite fields, although we shall mention some implications 
for Calabi-Yau manifolds defined over number fields. 

In this example, we see two types of mirror pairs. The first one is the 
maximally generic mirror pair $\{ X_{\Lambda}, W_{\lambda}\}$, where $X_{\Lambda}$ is 
the universal family of smooth projective Calabi-Yau hypersurfaces of 
degree $(n+1)$ in $\PP^{n}$ and $W_{\lambda}$ is the one parameter family 
of Calabi-Yau manifolds as constructed above. Note that $X_{\Lambda}$ and $Y_{\lambda}$ 
are parametrized by different parameter spaces (of different dimensions). 
The possible zeta symmetry in this case 
would then have to be a relation between certain generic property of the zeta function  
for $X\in X_{\Lambda}$ and the corresponding generic property of 
the zeta function for $W\in W_{\lambda}$. 

The second type of mirror pairs is the 
one parameter family of mirror pairs $\{ X_{\lambda}, W_{\lambda}\}$ parametrized 
by the same parameter $\lambda$. This is a stronger type of mirror pair than 
the first type. For $\lambda\in \C$, we say that $W_{\lambda}$ is a {\bf strong 
mirror} of $X_{\lambda}$. For such a strong mirror pair 
$\{ X_{\lambda}, W_{\lambda}\}$, we can really ask 
for the relation between the zeta function of $X_{\lambda}$ and the 
zeta function of $W_{\lambda}$. If $\lambda_1\not=\lambda_2$, $W_{\lambda_1}$ would 
not be called a strong mirror for $X_{\lambda_2}$, although they would be 
an usual {\bf weak mirror} pair. Apparently, we do not have a definition 
for a strong mirror pair in general, as there is not even a definition for 
a generic or weak mirror pair in general.

Let $\fq$ be a finite field of $q$ elements, where $q=p^r$ and $p$ is a prime. 
For a scheme $X$ of finite type of dimension $d$ over $\fq$, 
let $\# X(\fq)$ denote the number of $\fq$-rational 
points on $X$. Let 
$$Z(X,T) =\exp (\sum_{k=1}^{\infty} {T^k\over k} \# X(\fqk))\in 1+T\Z [[T]]$$
be the zeta function of $X$. It is well known that $Z(X,T)$ is a rational 
function in $T$ whose reciprocal zeros and reciprocal poles are 
Weil $q$-integers. Factor $Z(X,T)$ over the $p$-adic numbers $\C_p$ and write 
$$Z(X,T) =\prod_i (1-\alpha_iT)^{\pm 1}$$
in reduced form, where the algebraic integers $\alpha_i \in \C_p$. One knows that 
the slope ${\rm ord}_q(\alpha_i)$ is a rational number in the interval $[0,d]$.  
For two real numbers $s_1\leq s_2$, we define the slope $[s_1, s_2]$ part of 
$Z(X,T)$ to be the partial product 
\begin{equation}
Z_{[s_1,s_2]}(X,T) =\prod_{s_1\leq {\rm ord}_q(\alpha_i) \leq s_2} (1-\alpha_iT)^{\pm 1}.
\end{equation}
For a half open and half closed interval $[s_1, s_2)$, the 
slope $[s_1, s_2)$ part $Z_{[s_1, s_2)}(X,T)$ of $Z(X,T)$ is defined in a similar way. 
These are rational functions with coefficients in $\Z_p$ by the $p$-adic Weierstrass 
factorization. 
It is clear that we have the decomposition 
$$Z(X,T)= \prod_{i=0}^{d} Z_{[i, i+1)}(X, T).$$

Our main result of this paper is the following arithmetic mirror theorem. 

\begin{theorem} 
Assume that  $\lambda \in \fq$ such that 
$(X_{\lambda}, W_{\lambda})$ is a strong mirror pair of 
Calabi-Yau manifolds over $\fq$. 
For every positive integer $k$, we have the congruence formula 
$$\# X_{\lambda}(\fqk) \equiv  \# Y_{\lambda}(\fqk) \equiv \# W_{\lambda}(\fqk)   
 ~ ({\rm mod}~ q^k).$$
Equivalently, the slope $[0,1)$ part of the zeta function is the same 
for the mirror varieties $\{ X_{\lambda}, Y_{\lambda}, W_{\lambda}\}$: 
$$Z_{[0, 1)}(X_{\lambda},T) = Z_{[0,1)}(Y_{\lambda}, T)= Z_{[0,1)}(W_{\lambda}, T).$$ 
\end{theorem} 

We now discuss a few applications of this theorem. 
In terms of cohomology theory, this suggests that the semi-simplification of the 
DeRham-Witt cohomology ( in particular, the $p$-adic et\`ale cohomology) 
for $\{ X_{\lambda}, Y_{\lambda}, W_{\lambda}\}$ are all the same. A corollary of the 
above theorem is that the unit root parts (slope zero parts) of their zeta functions 
are the same: 
$$Z_{[0, 0]}(X_{\lambda},T) = Z_{[0,0]}(Y_{\lambda}, T)= Z_{[0,0]}(W_{\lambda}, T).$$
The $p$-adic variation of the rational function $Z_{[0, 0]}(X_{\lambda},T)$ 
as $\lambda$ varies is closely related to the mirror map 
which we do not discuss it here, but see \cite{Dw} for the case $n\leq 3$. 
From arithmetic point of view, the $p$-adic variation 
of the rational function $Z_{[0, 0]}(X_{\lambda},T)$ 
as $\lambda$ varies is explained 
by Dwork's unit root zeta function \cite{Dw2}. 
We briefly explain the connection here. 

Let $B$ be the parameter variety of $\lambda$ such that $(X_{\lambda}, W_{\lambda})$ 
form a strong mirror pair. Let $\Phi: X_{\lambda}\rightarrow B$ 
(resp. $\Psi: W_{\lambda}\rightarrow B$) be the projection to 
the base by sending $X_{\lambda}$ (resp. $W_{\lambda}$) to $\lambda$. 
The pair $(\Phi, \Psi)$ of morphisms to $B$ is called a {\bf strong mirror pair of 
morphisms to $B$}. Each of its fibres gives a strong mirror pair of 
Calabi-Yau manifolds.  
Recall that Dwork's unit root zeta function attached to the 
morphism $\Phi$ is defined to be the formal infinite product  
$$Z_{\rm unit}(\Phi, T) =\prod_{\lambda \in |B|}
Z_{[0,0]}(X_{\lambda}, T^{{\rm deg}(\lambda)}) \in 1+T\Z_p[[T]],$$
where $|B|$ denotes the set of closed points of $B$ over $\fq$. 
This unit root zeta function is no longer a rational function, 
but conjectured by Dwork in \cite{Dw2} and proved by the author 
in \cite{W4}\cite{W5}\cite{W6} to be a $p$-adic meromorphic function in $T$.  
The above theorem immediately implies 

\begin{corollary} Let $(\Phi, \Psi)$ be the above strong mirror pair of 
morphisms to the base $B$.  Then, their unit root zeta functions are the same:
$$Z_{\rm unit}(\Phi, T) =Z_{\rm unit}(\Psi, T).$$
\end{corollary}

If $\lambda$ is in a number field $K$, then Theorem 1.1 implies that the  
Hasse-Weil zeta functions of $X_{\lambda}$ and $Y_{\lambda}$ differ 
essentially by the L-function of a pure motive $M_n(\lambda)$ of 
weight $n-3$. That is, 
$$\zeta(X_{\lambda}, s) = \zeta(Y_{\lambda}, s) L(M_n(\lambda), s-1).$$
In the quintic case $n=4$, the pure weight $1$ motive $M_4(\lambda)$ would 
come from a curve. This curve has been constructed explicitly by 
Candelas, de la Ossa and Fernando-Rodriquez \cite{Ca}.
The relation between the Hasse-Weil zeta functions of 
$X_{\lambda}$ and $W_{\lambda}$ are similar, differing by a few more 
factors consisting of Tate twists of the Dedekind zeta function of $K$. 

Theorem 1.1 motivates the following more general conjecture. 

\begin{conjecture}[Congruence mirror conjecture]
Suppose that we are given a strong mirror pair $\{ X,Y\}$ of Calabi-Yau manifolds 
defined over $\fq$. Then, for every positive integer $k$, we have 
$$\# X(\fqk) \equiv  \# Y(\fqk)   
 ~ ({\rm mod}~ q^k).$$
Equivalently, 
$$Z_{[0, 1)}(X,T) = Z_{[0,1)}(Y, T).$$
Equivalently (by functional equation), 
$$Z_{(d-1, d]}(X,T) = Z_{(d-1,d]}(Y, T).$$  
\end{conjecture}

The condition in the congruence mirror 
conjecture is vague since one does not know at present 
an algebraic geometric definition of a strong mirror pair of 
Calabi-Yau manifolds, although one does know many examples 
such as the one given above. 
Thus, a major part of the problem is to make 
the definition of a strong mirror pair mathematically 
precise. For an additional evidence of the congruence mirror conjecture, 
see Theorem 6.2 which can be viewed as a generalization of Theorem 1.1. 
As indicated before, this conjecture implies that Dwork's unit root 
zeta functions for the two families forming a strong mirror pair are the same 
$p$-adic meromorphic functions. This means that under the strong 
mirror family involution, Dwork's unit root zeta function stays the same.

Just like the zeta function itself, its slope $[0,1)$ part 
$Z_{[0,1)}(X_{\lambda}, T)$ depends heavily on the algebraic 
parameter $\lambda$, not just on the topological properties of $X_{\lambda}$. 
This means that the congruence mirror conjecture is really a continuous type 
of arithmetic mirror symmetry. This continuous nature requires 
the use of a strong mirror pair, not just a generic mirror pair. 

Assume that $\{ X, Y\}$ forms a mirror pair, not necessarily 
a strong mirror pair. A different type of arithmetic mirror symmetry 
reflecting the Hodge symmetry, which is discrete and hence generic in nature,
is to look for a suitable quantum version $Z_Q(X,T)$ of 
the zeta function such that 
$$Z_Q(X,T)=Z_Q(Y,T)^{(-1)^d},$$
where $\{ X, Y\}$ is a mirror pair of Calabi-Yau manifolds 
over $\fq$ of dimension $d$. This relation cannot hold for the 
usual zeta function $Z(X,T)$ for obvious reasons, even for a strong mirror pair 
as it contradicts with the congruence mirror conjecture for odd $d$. No non-trivial 
candidate for $Z_Q(X,T)$ has been found. Here we propose a $p$-adic 
quantum version which would have the conjectural properties 
for most (and hence generic) mirror pairs. We will call our 
new zeta function to be the slope zeta function as it is based on 
the slopes of the 
zeros and poles.

\begin{definition}For a scheme $X$ of finite type over $\fq$, 
write as before 
$$Z(X,T) =\prod_i (1-\alpha_iT)^{\pm 1}$$
in reduced form, where $\alpha_i \in \C_p$.
Define the slope zeta function of $X$ to be the 
two variable function 
\begin{equation}
S_p(X, u, T) =\prod_i (1-u^{{\rm ord}_q(\alpha_i)}T)^{\pm 1}.
\end{equation}
\end{definition}

Note that 
$$\alpha_i  = q^{{\rm ord}_q(\alpha_i)} \beta_i,$$
where $\beta_i$ is a $p$-adic unit. Thus, the slope zeta function 
$S_p(X, u, T)$ is obtained from the $p$-adic factorization of $Z(X,T)$
by dropping the $p$-adic unit parts of the roots and replacing $q$ by the variable 
$u$. This is not always a rational function in $u$ and $T$. 
It is rational if all slopes are integers. 
Note that the definition of the slope zeta function is independent of 
the choice of the ground field $\fq$ where $X$ is defined. It depends 
only on $X\otimes \fqbar$ and thus is also a geometric invariant. 
It would be interesting to see if there is a 
diophantine interpretation of the slope zeta function.

If $X$ is a scheme of finite type over $\Z$, then for each prime number 
$p$, the reduction $X\otimes \fp$ has the $p$-adic slope zeta function 
$S_p(X\otimes \fp, u, T)$. At the first glance, one might think that this 
gives infinitely many discrete invariants for $X$ 
as the set of prime numbers is infinite.  
However, it can be shown that the set $\{ S_p(X\otimes \fp, u, T) | p \ {\rm prime}\}$ 
contains only finitely many distinct elements. 
In general, it is a very interesting but difficult problem to determine this set 
$\{ S_p(X\otimes \fp, u, T) | p \ {\rm prime}\}$.

Suppose that $X$ and $Y$ form a mirror pair of $d$-dimensional Calabi-Yau manifolds 
over $\fq$. For simplicity and for comparison with the Hodge theory, 
we always assume in this paper that $X$ and $Y$ can be lifted to characteristic zero 
(to the Witt ring of $\fq$). In this good reduction case, 
the modulo $p$ Hodge numbers equal the characteristic zero Hodge numbers. 
Taking $u=1$ in the definition of the slope zeta function, 
we see that the specialization $S_p(X, 1, T)$ 
already satisfies the desired relation
$$S_p(X, 1, T) =(1-T)^{-e(X)}=(1-T)^{-(-1)^de(Y)} =S_p(Y, 1, T)^{(-1)^d}.$$
This suggests that there is a chance 
that the slope zeta function might satisfy the desired slope mirror  symmetry
\begin{equation}
S_p(X,u, T)={S_p(Y,u,T)^{(-1)^d}}.
\end{equation}
In section $7$, we shall show that the slope zeta function satisfies a functional 
equation. Furthermore, the expected slope mirror symmetry does hold 
if both $X$ and $Y$ are ordinary. If either $X$ or $Y$ is not ordinary, 
the expected slope mirror symmetry 
is unlikely to hold in general.

If $d\leq 2$, the congruence mirror conjecture implies that 
the slope zeta function does satisfy the expected slope mirror symmetry 
for a strong mirror pair $\{ X, Y\}$,  
whether $X$ and $Y$ are ordinary 
or not. For $d\geq 3$, we believe that the slope zeta function
is still a little bit too strong for the expected 
symmetry to hold in general, even if $\{ X, Y\}$ forms a strong mirror pair. 
And it should not be too hard to find 
a counter-example although we have not done so. However, we believe 
that the expected slope mirror symmetry holds for a sufficiently generic 
pair of $3$-dimensional Calabi-Yau manifolds. 

\begin{conjecture}[Slope mirror conjecture] 
Suppose that we are given a maximally generic mirror pair
$\{ X, Y\}$ of $3$-dimensional 
Calabi-Yau manifolds defined over $\fq$. Then, we have the slope mirror symmetry 
\begin{equation}
S_p(X,u, T)={1\over S_p(Y,u,T)}
\end{equation}
for generic $X$ and generic $Y$. 
\end{conjecture}

A main point of this conjecture is that it holds for all prime numbers $p$. 
For arbitrary $d\geq 4$, the corresponding slope mirror conjecture 
might be false for some prime numbers $p$, 
but it should be true for all primes $p\equiv 1 ~({\rm mod}~D)$ for some 
positive integer $D$ depending on the mirror family, if the family comes from 
the reduction modulo $p$ of a family defined over a number field. In the case $d\leq 3$, 
one could take $D=1$ and hence get the above conjecture. 

Again the condition in the slope mirror conjecture 
is vague as it is not presently known 
an algebraic geometric definition of a mirror family, 
although many examples are known in the toric setting. 
In a future paper, using the results in \cite{W1}\cite{W2}, 
we shall prove that the slope mirror conjecture 
holds in the toric hypersurface case if $d\leq 3$. For example, 
if $X$ is a generic quintic hypersurface, then $X$ is ordinary 
by the results in \cite{IL}\cite{W1} for every $p$ and thus one finds 
$$S_p(X\otimes \fp, u,T) ={(1-T)(1-uT)^{101}(1-u^2T)^{101}(1-u^3T) 
\over (1-T)(1-uT)(1-u^2T)(1-u^3T)}.$$
This is independent of $p$. 
Note that we do not know if the one parameter subfamily $X_{\lambda}$ 
is generically ordinary for every $p$. The ordinary property for every $p$ was 
established only for the universal family of hypersurfaces, not 
for a one parameter subfamily of hypersurfaces such as $X_{\lambda}$.  
If $Y$ denotes the generic mirror of $X$, then by the results in \cite{W1} 
\cite{W2}, $Y$ is ordinary for every $p$ and thus  
we obtain
$$S_p(Y\otimes \fp, u,T) ={ (1-T)(1-uT)(1-u^2T)(1-u^3T)\over 
(1-T)(1-uT)^{101}(1-u^2T)^{101}(1-u^3T)}.$$
Again, it is independent of $p$. 
The slope mirror conjecture holds in this example.

For a mirror pair over a number field, we have the 
following harder conjecture.

\begin{conjecture}[Slope mirror conjecture over $\Z$] 
Let $\{ X, Y\}$ be two schemes of finite type over $\Z$ 
such that their generic fibres $\{ X\otimes \Q, Y\otimes \Q\}$ 
form a usual (weak) mirror pair of $d$-dimensional 
Calabi-Yau manifolds defined over $\Q$. Then there are infinitely many 
prime numbers $p$ (with positive density) such that  
$$S_p(X\otimes \fp, u, T)=S_p(Y\otimes \fp, u,T)^{(-1)^d}.$$
\end{conjecture}

{\bf Remarks}. If one uses the weight $2\log_q|\alpha_i|$ instead of the slope 
${\rm ord}_q\alpha_i$, 
where $|\cdot |$ denotes the complex absolute value, one can define a two 
variable weight zeta function 
in a similar way. It is easy to see that the resulting weight zeta function 
does not satisfy the desired symmetry as the weight has nothing to do with 
the Hodge symmetry, while the slopes are related to the Hodge numbers as 
the Newton polygon (slope polygon) lies above the Hodge polygon. 

In practice, one is often given a mirror pair of singular Calabi-Yau orbifolds, 
where there may not exist a smooth crepant resolution. In such a case, one could  
define an orbifold zeta function, which would be equal to the zeta function of 
the smooth crepant resolution whenever such a resolution exists. 
Similar results and conjectures should carry over to such orbifold zeta 
functions.

\section{A counting formula via Gauss sums}

Let $V_1,\cdots, V_m$ be $m$ distinct lattice points in $\Z^n$. 
For $V_j=(V_{1j}, \cdots, V_{nj})$, write 
$$x^{V_j} = x_1^{V_{1j}}\cdots x_n^{V_{nj}}.$$
Let $f$ be the Laurent polynomial in $n$ variables written in the form: 
\[
f(x_1,\cdots, x_n)=\sum_{j=1}^m{a_j}x^{V_j}, { a_j}\in \fq, 
\]
where not all $a_j$ are zero. 
Let $M$ be the $n\times m$ matrix 
\[
M=(V_1, \cdots, V_m), 
\]
where 
each $V_j$ is written as a column vector. 
Let $N_f^*$ denote the number of $\fq$-rational points 
on the affine toric hypersurface $f=0$ in $\G_m^n$. 
If each $V_j\in \Z_{\geq 0}^n$, we let $N_f$ denote the 
number of $\fq$-rational points on the affine hypersurface 
$f=0$ in $\A^n$. We first derive a well known formula 
for both $N_f^*$ and $N_f$ in terms of Gauss sums.

For this purpose, we now recall the definition of Gauss sums. 
Let $\fq$ be the finite field of $q$ elements, where $q=p^r$ and $p$ 
is the characteristic of $\fq$. 
Let $\chi$ be the Teichm\"uller character of the multiplicative group $\fq^{*}$. 
For $a\in \fq^*$, the value $\chi(a)$ is just the $(q-1)$-th root of unity 
in the $p$-adic field ${\C}_p$ 
such that $\chi(a)$ modulo $p$ reduces to $a$.  Define the $(q-2)$ 
Gauss sums over  $\fq$ by
$$G(k)=\sum_{a\in \fq^*}\chi(a)^{-k}\zeta_p^{{\rm Tr}(a)}\ \ (1\leq k\leq q-2), $$
where $\zeta_p$ is a primitive $p$-th root of unity in $\C_p$ and 
${\rm Tr}$ denotes the trace map from $\fq$ to the prime field $\fp$. 

\begin{lemma} For all $a\in \fq$, the Gauss sums satisfy the following 
interpolation relation 
$$\zeta_p^{ {\rm Tr}(a)} =\sum_{k=0}^{q-1} {G(k) \over q-1} \chi(a)^k,$$
where 
$$G(0)=q-1, \ G(q-1)=-q.$$
\end{lemma}

{\bf Proof}. By the Vandermonde determinant, there are numbers $C(k)$ ($0\leq k\leq q-1$) 
such that for all $a\in \fq$, one has 
$$\zeta_p^{ {\rm Tr}(a)} =\sum_{k=0}^{q-1} {C(k) \over q-1} \chi(a)^k.$$
It suffices to prove that $C(k)=G(k)$ for all $k$. 
Take $a=0$, one finds that $C(0)/(q-1)=1$. This proves that 
$C(0)=q-1=G(0)$. For $1 \leq k\leq q-2$, one computes that 
$$G(k)=\sum_{a\in \fq^*}\chi(a)^{-k}\zeta_p^{{\rm Tr}(a)} ={C(k)\over q-1}(q-1) =C(k).$$
Finally, 
$$0 = \sum_{a\in \fq}\zeta_p^{{\rm Tr}(a)} ={C(0)\over q-1}q + {C(q-1)\over q-1}(q-1).$$
This gives $C(q-1)=-q =G(q-1)$. The lemma is proved. 
\medskip 

We also need to use the following classical theorem of Stickelberger. 

\begin{lemma}Let $0\leq k\leq q-1$. Write 
$$k=k_0+k_1p+\cdots +k_{r-1}p^{r-1}$$
in $p$-adic expansion, where $0\leq k_i \leq p-1$. 
Let $\sigma(k)=k_0+\cdots +k_{r-1}$ be the sum of the 
$p$-digits of $k$. Then, 
$${\rm ord}_p G(k) ={\sigma(k)\over p-1}.$$
\end{lemma}

Now we turn to deriving a counting formula for $N_f$ in terms of Gauss sums. 
Write $W_j =(1, V_j)\in {\Z}^{n+1}$. Then, 
$$x_0f =\sum_{j=1}^m a_j x^{W_j}=\sum_{j=1}^m a_j x_0x_1^{V_{1j}}\cdots x_n^{V_{nj}},$$
where $x$ now has $n+1$ variables $\{ x_0, \cdots, x_n\}$. 
Using the formula 
$$\sum_{t\in \fq} t^k =\cases{0, & if $(q-1)\not| k$, \cr 
                              q-1, & if $(q-1)|k$ and $k>0$, \cr 
                              q, & if $k=0$, \cr}$$
one then calculates that 
\begin{eqnarray}
qN_f&=&\sum_{x_0, \cdots, x_n \in \fq}\zeta_p^{{\rm Tr}(x_0f(x))} \cr
&=&\sum_{x_0, \cdots, x_n \in \fq}\prod_{j=1}^m
\zeta_p^{ {\rm Tr}(a_j x^{W_j})} \cr
&=&\sum_{x_0, \cdots, x_n\in \fq}\prod_{j=1}^m
\sum_{k_j=0}^{q-1}{G(k_j)\over q-1} \chi(a_j)^{k_j}\chi(x^{W_j})^{k_j}\cr
&=&\sum_{k_1=0}^{q-1}\cdots \sum_{k_m=0}^{q-1}(\prod_{j=1}^m{G(k_j)\over q-1} 
\chi(a_j)^{k_j})
\sum_{x_0, \cdots, x_n \in \fq}\chi(x^{k_1W_1 +\cdots +k_mW_m})
\cr
&=&\sum_{\sum_{j=1}^m k_jW_j\equiv 0({\rm mod}~q-1)}{(q-1)^{s(k)}q^{n+1-s(k)} \over (q-1)^m}
\prod_{j=1}^m\chi(a_j)^{k_j} G(k_j), 
\end{eqnarray}
where $s(k)$ denotes the number of non-zero entries in $k_1W_1+\cdots +k_mW_m$. 

Similarly, one calculates that 
\begin{eqnarray}
qN_f^*&=&\sum_{x_0\in \fq, x_1, \cdots, x_n \in \fq^*}\zeta_p^{{\rm Tr}(x_0f(x))} \cr
&=&(q-1)^n + \sum_{x_0, \cdots, x_n \in \fq^*}\prod_{j=1}^m
\zeta_p^{ {\rm Tr}(a_j x^{W_j})} \cr
&=&(q-1)^n +\sum_{\sum_{j=1}^m k_jW_j\equiv 0({\rm mod}~q-1)}{(q-1)^{n+1} \over (q-1)^m}
\prod_{j=1}^m\chi(a_j)^{k_j} G(k_j). 
\end{eqnarray}
We shall use these two formulas to study the number of $\fq$-rational points 
on certain hypersurfaces in next two sections.

\section{Rational points on Calabi-Yau hypersurfaces}

In this section, we apply formula (7) to compute the number of 
$\fq$-rational points on the projective hypersurface $X_{\lambda}$ in $\PP^n$ defined 
by 
$$f(x_1,\cdots, x_{n+1})= x_1^{n+1}+\cdots +x_{n+1}^{n+1}+\lambda x_1\cdots x_{n+1}=0,$$
where $\lambda$ is an element of $\fq^*$. We shall handle the easier 
case $\lambda=0$ separately. 
Let $M$ be the $(n+2)\times (n+2)$ matrix 
\begin{equation}
M= \pmatrix{1 &1 &1 &\cdots &1 &1  \cr
            n+1 &0 &0 &\cdots &0 &1  \cr
             0 &n+1 &0 &\cdots &0 &1  \cr
             \vdots & \vdots &\vdots & \cdots & \vdots & \vdots \cr
             0&0 &0 &\cdots  &n+1 &1} 
\end{equation}
Let $k=(k_1,\cdots, k_{n+2})$ written as a column vector. 
Let $N_f$ denote the number of $\fq$-rational points on the 
affine hypersurface $f=0$ in ${\A}^{n+1}$. 
By formula (7), we deduce that 
$$qN_f =\sum_{Mk\equiv 0 ({\rm mod}~ q-1)} {(q-1)^{s(k)}q^{n+2-s(k)}\over (q-1)^{n+2}} 
(\prod_{j=1}^{n+2}G(k_j))\chi(\lambda)^{k_{n+2}},$$
where $s(k)$ denotes the number of non-zero entries in $Mk\in {\Z}^{n+2}$. 
The number of $\fq$-rational points on the projective hypersurface $X_{\lambda}$ is then 
given by the formula
$${N_f-1 \over q-1} = {-1\over q-1} + 
\sum_{Mk\equiv 0 ({\rm mod}~ q-1)} {q^{n+1-s(k)}\over (q-1)^{n+3-s(k)}} 
(\prod_{j=1}^{n+2}G(k_j))\chi(\lambda)^{k_{n+2}}.$$  
If $k=(0, \cdots, 0, q-1)$, then $Mk=(q-1,\cdots, q-1)$ and $s(k)=n+2$. 
In this case, the corresponding term in the above expression is 
$-(q-1)^n$ which is $(-1)^{n-1}$ modulo $q$. 
If $k=(0,..., 0)$, then $s(k)=0$ and the corresponding term is 
$q^{n+1}/(q-1)$ which is zero modulo $q$. 

Thus, we obtain the congruence formula modulo $q$: 
$${N_f-1 \over q-1} \equiv 1 +(-1)^{n-1} + 
{\mathop{{\sum}^*}_{Mk\equiv 0 ({\rm mod}~ q-1)}} {q^{n+1-s(k)}\over (q-1)^{n+3-s(k)}} 
(\prod_{j=1}^{n+2}G(k_j))\chi(\lambda)^{k_{n+2}},$$
where $\sum^*$ means summing over all those solutions $k=(k_1,\cdots, k_{n+2})$ 
with $0\leq k_i\leq q-1$, 
$k\not=(0,\cdots, 0)$, and $k\not=(0,\cdots, 0, q-1)$. 

\begin{lemma} If $k\not=(0,\cdots, 0)$, then $\prod_{j=1}^{n+2}G(k_j)$ 
is divisible by $q$.
\end{lemma}

{\bf Proof}. Let $k$ be a solution of $Mk\equiv 0({\rm mod}~q-1)$ such that 
$k\not=(0,\cdots, 0)$. Then, there are positive integers $\ell_0, \cdots, \ell_{r-1}$ 
such that 
$$k_1+\cdots +k_{n+2} =(q-1)\ell_0,$$
$$<pk_1>+\cdots +<pk_{n+2}> =(q-1)\ell_1,$$
$$\cdots$$
$$<p^{r-1}k_1>+\cdots +<p^{r-1}k_{n+2}> =(q-1)\ell_{r-1},$$
where $<pk_1>$ denotes the unique integer in $[0, q-1]$ 
congruent to $pk_1$ modulo $(q-1)$ and which is $0$ (resp. $q-1$) 
if $pk_1=0$ (resp., if $pk_1$ is a positive multiple of $q-1$). 
By the Stickelberger theorem, we deduce that 
$${\rm ord}_p\prod_{j=1}^{n+2}G(k_j) ={\sum_j \sigma(k_j) \over p-1} 
={1\over q-1}\sum_{i=0}^{r-1}(q-1)\ell_i =\sum_{i=0}^{r-1}\ell_i.$$
Since $\ell_i\geq 1$, it follows that 
$${\rm ord}_q\prod_{j=1}^{n+2}G(k_j) ={1\over r} \sum_{i=0}^{r-1}\ell_i \geq 1$$
with equality holding if and only if all $\ell_i=1$. 
The lemma is proved.

Using this lemma and the previous congruence formula, we deduce 

\begin{lemma} Let $\lambda\in \fq^*$. We have the congruence formula modulo $q$: 
$$\# X_{\lambda}(\fq) \equiv 1 +(-1)^{n-1} + 
{\mathop{{\sum}^*}_{Mk\equiv 0 ({\rm mod}~ q-1) \atop s(k)=n+2}} {1\over q(q-1)} 
(\prod_{j=1}^{n+2}G(k_j))\chi(\lambda)^{k_{n+2}}.$$
\end{lemma}

\section{Rational points on the mirror hypersurfaces}

In this section, we apply formula (8) to compute the number of 
$\fq$-rational points on the affine toric hypersurface in $\G_m^n$ defined 
by the Laurent polynomial equation 
$$g(x_1,\cdots, x_{n})= x_1+\cdots +x_n+{1\over x_1\cdots x_{n}}+\lambda=0,$$
where $\lambda$ is an element of $\fq^*$. 
Let $N$ be the $(n+1)\times (n+2)$ matrix 
\begin{equation}
N= \pmatrix{1 &1 &\cdots &1 &1 &1 \cr
            1 &0 &\cdots &0 &-1&0  \cr
             0 &1 &\cdots &0 &-1&0  \cr
             \vdots & \vdots & \cdots & \vdots & \vdots \vdots \cr
             0&0  &\cdots  &1 &-1 &0} 
\end{equation}
Let $k=(k_1,\cdots, k_{n+2})$ written as a column vector. 
By formula (8), we deduce that 
$$qN_g^* =(q-1)^n+\sum_{Nk\equiv 0 ({\rm mod}~ q-1)} {1 \over (q-1)} 
(\prod_{j=1}^{n+2}G(k_j))\chi(\lambda)^{k_{n+2}},$$
where $k=(k_1,\cdots, k_{n+2})$ with $0\leq k_i\leq q-1$.

The contribution of those trivial terms $k$ (where each $k_i$ is either 
$0$ or $q-1$) is given by 
$${1\over q-1}\sum_{s=0}^{n+2}(-q)^s (q-1)^{n+2-s}{n+2\choose s} ={(-1)^n\over q-1}.$$
Since 
$$(q-1)^n +{(-1)^n\over q-1} = {(q-1)^{n+1}+(-1)^n \over q-1} \equiv q(n+1)(-1)^{n-1} 
({\rm mod} q^2),$$
we deduce 

\begin{lemma} For $\lambda\in \fq^*$, 
we have the following congruence formula modulo $q$:
$$N_g^* \equiv (n+1)(-1)^{n-1} +{\mathop{{\sum}'}_{Nk\equiv 0 ({\rm mod}~ q-1)}} 
{1 \over q(q-1)} 
(\prod_{j=1}^{n+2}G(k_j))\chi(\lambda)^{k_{n+2}},$$
where $\sum'$ means summing over all those non-trivial solutions $k$. 
\end{lemma}

\section{The mirror congruence formula}

\begin{theorem} For $\lambda\in \fq^*$, we have the congruence formula 
$$\# X_{\lambda}(\fq) \equiv N_g^* +1-n(-1)^{n-1} ~({\rm mod}~ q).$$
\end{theorem}

{\bf Proof}. If $k$ is a non-trivial solution of $Nk \equiv 0 ({\rm mod} q-1)$, 
then we have 
$$k_1\equiv k_2 \equiv \cdots \equiv k_n \equiv k_{n+1} ({\rm mod} q-1)$$ 
and 
$$k_1+\cdots +k_{n+1}+k_{n+2} \equiv 0 ({\rm mod} q-1).$$
Since $k$ is non-trivial, we must have 
$$0<k_1=k_2=\cdots =k_{n+1}<q-1,$$
$$k_1+\cdots +k_{n+2} =(n+1)k_1+k_{n+2}=(n+1)k_2+k_{n+2}=\cdots \equiv 0 ({\rm mod} q-1).$$
This gives all solutions of the equation $Mk\equiv 0 ({\rm mod}q-1)$ 
with $k_1=\cdots =k_{n+1}$, $0<k_1<q-1$ and $s(k)=n+2$. 
The corresponding terms for these $k$'s in $(N_f-1)/(q-1)$ and $N_g^*$ 
are exactly the same. 

A solution of $Mk\equiv 0 ({\rm mod}~q-1)$ is called {\bf admissible} if 
$s(k)=n+2$ and its first $k+1$ coordinates $\{k_1,\cdots, k_{n+1}\}$ contain 
at least two distinct elements.  
The above results show that we have  
$${N_f -1 \over q-1} -1 -(-1)^{n-1} -(N_g^* -(n+1)(-1)^{n-1})$$
$$\equiv \sum_{{\rm admissible}~ k} 
{1 \over q(q-1)} 
(\prod_{j=1}^{n+2}G(k_j))\chi(\lambda)^{k_{n+2}} ~({\rm mod}~ q).$$
This congruence together with the following lemma 
completes the proof of the theorem.

\begin{lemma}If $k$ is an admissible solution of $Mk\equiv 0 ({\rm mod}~q-1)$, 
then 
$${\rm ord}_q (\prod_{j=1}^{n+2}G(k_j)) \geq 2.$$
\end{lemma}

{\bf Proof}. If $k$ is an admissible solution, then $<pk>,\cdots, <p^{r-1}k>$ are
also admissible solutions. For each $1\leq i\leq n+1$, write 
$$(n+1)k_i +k_{n+2} =(q-1)\ell_i,$$
where $\ell_i$ is a positive integer. 
Adding these equations together, we get 
$$(n+1)(k_1+\cdots +k_{n+1}) +(n+1)k_{n+2}=(q-1)(\ell_1+\cdots +\ell_{n+1}).$$
Thus, the integer 
$${k_1+\cdots +k_{n+2}\over q-1} ={\ell_1+\cdots +\ell_{n+1} \over n+1} 
=\ell \in {\bf Z}_{>0}.$$
It is clear that $\ell =1$ if and only if each $\ell_i=1$ which would imply that 
$k_1=\cdots =k_{n+1}$ contradicting with the admissibility of $k$. 
Thus, we must have that $\ell \geq 2$. 

Similarly, for each $0\leq i\leq r-1$, we have 
$$<p^ik_1> +\cdots +<p^ik_{n+2}> =(q-1)j_i,$$
where $j_i\geq 2$ is a positive integer. 
We conclude that
$${\rm ord}_q (\prod_{j=1}^{n+2}G(k_j)) ={j_0+\cdots +j_{r-1} \over r} \geq 2.$$
The lemma is proved. 

\section{Rational points on the projective mirror}

Let ${\Delta}$ be the convex integral polytope associated with the 
Laurent polynomial $g$. It is the $n$-dimensional simplex in $\R^n$ 
with the following vertices:
$$\{ e_1, \cdots, e_n, -(e_1+\cdots +e_n)\}, $$
where the $e_i$'s are the standard unit vectors in $\R^n$. 

Let $\PP_{\Delta}$ be the projective toric variety 
associated with the polytope $\Delta$, which contains $\G_m^n$ 
as an open dense subset. Let $Y_{\lambda}$ be the projective closure 
in $\PP_{\Delta}$ of the affine toric hypersurface  $g=0$ in $\G_m^n$. 
The variety $Y_{\lambda}$ is then a projective toric hypersurface 
in $\PP_{\Delta}$. We are interested in the number of $\fq$-rational 
points on $Y_{\lambda}$. 

The toric variety $\PP_{\Delta}$ has the following disjoint decomposition:
$$\PP_{\Delta} = \bigcup_{\tau \in \Delta} \PP_{\Delta, \tau},$$
where $\tau$ runs over all non-empty faces of $\Delta$ and each 
$\PP_{\Delta, \tau}$ is isomorphic to the torus $\G_m^{{\rm dim}\tau}$.
Accordingly, the projective toric hypersurface $Y_{\lambda}$ has the 
corresponding disjoint decomposition 
$$Y_{\lambda} = \bigcup_{\tau \in \Delta} Y_{\lambda, \tau}, 
\ Y_{\lambda, \tau} = Y_{\lambda}\cap {\PP}_{\Delta, \tau}.$$
For $\tau=\Delta$, the subvariety $Y_{\lambda, \Delta}$ is simply 
the affine toric hypersurface defined by $g=0$ in $\G_m^n$. 
For zero-dimensional $\tau$, $Y_{\lambda, \tau}$ is empty. 
For a face $\tau$ with $1\leq {\rm dim}\tau \leq n-1$, 
one checks that $Y_{\lambda, \tau}$ is isomorphic to the 
affine toric hypersurface in $\G_m^{{\rm dim}\tau}$ defined by 
$$1+x_1+\cdots +x_{{\rm dim}\tau}=0.$$
For such a $\tau$, the inclusion-exclusion principle shows that 
$$\# Y_{\lambda, \tau}(\fq) = q^{{\rm dim}\tau-1} -
{{\rm dim}\tau \choose 1} q^{{\rm dim}\tau-2} +\cdots +
(-1)^{{\rm dim}\tau-1}{{\rm dim}\tau \choose {\rm dim}\tau-1}.$$
Thus, 
$$\# Y_{\lambda, \tau}(\fq) ={1\over q}((q-1)^{{\rm dim}\tau} +
(-1)^{{\rm dim}\tau +1}).$$
This formula holds even for zero-dimensional $\tau$ as both sides 
would then be zero. 

Putting these calculations together, we deduce that 
$$\# Y_{\lambda}(\fq) = N_g^* -{(q-1)^n +(-1)^{n+1}\over q} 
+\sum_{\tau \in \Delta}{1\over q}((q-1)^{{\rm dim}\tau} +
(-1)^{{\rm dim}\tau +1}),$$
where $\tau$ runs over all non-empty faces of $\Delta$ 
including $\Delta$ itself.  
Since $\Delta$ is a simplex, one computes that
$$\sum_{\tau \in \Delta}((q-1)^{{\rm dim}\tau} +
(-1)^{{\rm dim}\tau +1}) ={q^{n+1}-1\over q-1} +(-1)={q(q^n-1)\over q-1}.$$
This implies that  
\begin{equation}
\# Y_{\lambda}(\fq) = N_g^* -{(q-1)^n +(-1)^{n+1}\over q} 
+{q^n-1\over q-1}.
\end{equation} 
This equality holds for all $\lambda \in \fq$, including the case $\lambda=0$. 
Reducing modulo $q$, we get 
\begin{equation}
\# Y_{\lambda}(\fq) \equiv N_g^* +1 -n(-1)^{n-1} ~({\rm mod}~ q). 
\end{equation}
This and Theorem 5.1 prove the case $\lambda \not=0$ of 
the following theorem.

\begin{theorem} For every finite field $\fq$ with $\lambda \in \fq$, 
we have the congruence formula 
$$\# X_{\lambda}(\fq)   \equiv \# Y_{\lambda}(\fq)~({\rm mod}~ q). $$ 
If furthermore, $\lambda\in \fq$ such that $g$ is $\Delta$-regular and $W_{\lambda}$ 
is a mirror manifold of $X_{\lambda}$, 
then 
$$\# Y_{\lambda}(\fq)   \equiv \# W_{\lambda}(\fq)~({\rm mod}~ q). $$  
\end{theorem}

{\bf Proof}. For the first part, it remains to check the case $\lambda=0$. The proof is 
similar and in fact somewhat simpler than the case $\lambda\not=0$. 
We give an outline. Since $\lambda=0$, we can take $k_{n+2}=0$ 
in the calculations of $N_f$ and $N_g^*$. One finds then 
$$\# X_0(\fq) \equiv 1 + 
{\mathop{{\sum}^*}_{Mk\equiv 0 ({\rm mod}~ q-1) \atop s(k)=n+2}} {1\over q(q-1)} 
(\prod_{j=1}^{n+2}G(k_j)),$$
where  $\sum^*$ means summing over all those solutions $k=(k_1,\cdots, k_{n+1},0)$ 
with $0\leq k_i\leq q-1$ and $k\not=(0,\cdots, 0)$. 

Similarly, one computes that 
$$N_g^* \equiv n(-1)^{n-1} +{\mathop{{\sum}'}_{Nk\equiv 0 ({\rm mod}~ q-1)}} 
{1 \over q(q-1)} 
(\prod_{j=1}^{n+2}G(k_j)),$$
where $\sum'$ means summing over all those non-trivial solutions $k$ 
with $k_{n+2}=0$.  
By (12), we deduce  
$$\# Y_0(\fq) \equiv 1 +{\mathop{{\sum}'}_{Nk\equiv 0 ({\rm mod}~ q-1)}} 
{1 \over q(q-1)} 
(\prod_{j=1}^{n+2}G(k_j)).$$
As before, one checks that 
$${\mathop{{\sum}^*}_{Mk\equiv 0 ({\rm mod}~ q-1) \atop s(k)=n+2}} 
(\prod_{j=1}^{n+2}G(k_j)) \equiv {\mathop{{\sum}'}_{Nk\equiv 0 ({\rm mod}~ q-1)}}
(\prod_{j=1}^{n+2}G(k_j)) ({\rm mod}~q^2).$$
The first part of the theorem follows. 

To prove the second part of the theorem, let $\Delta^*$ be the dual polytope of $\Delta$. 
One checks that $\Delta^*$ is the simplex in $\R^n$ with the vertices 
$$(n+1)e_i-\sum_{j=1}^n e_j ~(i=1,..., n), ~ -\sum_{j=1}^n e_j.$$
This is the $(n+1)$-multiple of a basic (regular) simplex in $\R^n$. In particular, 
the codimension $1$ faces of $\Delta^*$ are $(n+1)$-multiples of a basic simplex 
in $\R^{n-1}$. 
By the parrallel hyperplane decomposition in \cite{KK}, one deduces that the codimension 
$1$ faces of $\Delta^*$ have a triangulation into basic simplices. 
Fix such a triangulation  
which produces a smooth crepant resolution $\phi: W_{\lambda} \rightarrow Y_{\lambda}$. 
One checks \cite{Ba2} that for 
each point $y\in Y_{\lambda}(\fq)$, 
the fibre $\phi^{-1}(\lambda)$ is stratified by affine spaces over $\fq$. Since the 
fibres are connected, it follows that 
the number of $\fq$-rational points on $\phi^{-1}(\lambda)$ is congruent to $1$ 
modulo $q$. Thus, modulo $q$, we have the congruence 
$$\# W_{\lambda}(\fq) \equiv \sum_{y\in Y_{\lambda}(\fq)} \phi^{-1}(\lambda)(\fq) 
\equiv \sum_{y\in Y_{\lambda}(\fq)} 1 = \# Y_{\lambda}(\fq).$$
The proof is complete.

In terms of zeta functions, the above theorem says that the slope $[0,1)$ part 
of the zeta function for $X_{\lambda}$ equals the slope $[0,1)$ part of the 
zeta function for $Y_{\lambda}$. 

The above elementary calculations can be used to treat some other examples 
of toric hypersurfaces and complete intersections. In a forthcoming joint work with 
Lei Fu, we can prove the following generalization. 

\begin{theorem}Let $X$ be a smooth connected Calabi-Yau variety 
defined over the ring $W(\fq)$ of Witt vectors of $\fq$. Let $G$ 
be a finite group acting on $X$. Assume that $G$ fixes the non-zero 
global section of the canonical bundle of $X$. Then, for each positive integer $k$, 
we have the 
congruence formula 
$$\# (X\otimes\fq)(\fqk) \equiv \# (X/G \otimes \fq)(\fqk) ({\rm mod}~q^k).$$ 
\end{theorem}

Strictly speaking, this is not a complete  generalization of Theorem 6.1 yet, 
since Theorem 6.1 includes singular cases as well.

\section{Applications to zeta functions}

In this section, we compare the two zeta functions $Z(X_{\lambda}, T)$ and 
$Z(Y_{\lambda}, T)$, where $\{ X_{\lambda}, Y_{\lambda}\}$ is our strong mirror pair.  

First, we recall what is known about $Z(X_{\lambda}, T)$.  
Let $\lambda\in \fq$ such that $X_{\lambda}$ is smooth projective. 
By the Weil conjectures, the zeta function of $X_{\lambda}$ over $\fq$ has the 
following form
\begin{equation}
Z(X_{\lambda}, T) = {{P(\lambda, T)^{(-1)^n}}\over (1-T)(1-qT)\cdots (1-q^{n-1}T)},
\end{equation}
where $P(\lambda, T)\in 1+T\Z [T]$ is a polynomial of degree $n(n^n-(-1)^n)/(n+1)$, 
pure of weight $n-1$. By the results in \cite{IL}\cite{W1}, the universal family of 
hypersurfaces of degree $n+1$ is generically ordinary for every $p$ (Mazur's conjecture). 
However, we do not know if the one parameter family $X_{\lambda}$ of 
hypersurfaces is generically ordinary for every $p$. Thus, we raise 

\begin{question}Is the one parameter family $X_{\lambda}$ of degree $n+1$ 
hypersurfaces in $\PP^n$ generically ordinary for every prime number $p$ 
not dividing $(n+1)$? 
\end{question}

The answer is yes if $p \equiv 1~({\rm mod}~n+1)$ since the fibre for $\lambda=0$ 
is already ordinary if $p \equiv 1~({\rm mod}~n+1)$. It is also true if $n\leq 3$. 
The first unknown case is when $n=4$, the quintic case.  

Next, we recall what is known about $Z(Y_{\lambda}, T)$.   
Let $\lambda\in \fq$ such that $g$ is $\Delta$-regular. 
This is equivalent to assuming that $\lambda^n \not=(n+1)^{n+1}$. 
Then, the zeta function of the affine toric hypersurface $g=0$ over $\fq$ in $\G_m^n$  
has the following form (see \cite{W3}) 
$$Z(g, T) = {{Q(\lambda, T)^{(-1)^n}}} 
\prod_{i=0}^{n-1}(1-q^iT)^{(-1)^{n-i}{n\choose i+1}},$$
where $Q(\lambda, T)\in 1+T\Z [T]$ is a polynomial of degree $n$, pure of weight $n-1$.
The product of the trivial factors in $Z(g,T)$ is simply the zeta function 
of this sequence 
$${(q^k-1)^n +(-1)^{n+1} \over q^k}, \ k=1,2,\cdots.$$  
From this and (11), 
one deduces that the zeta function of the projective toric hypersurface 
$Y_{\lambda}$ has the form 
\begin{equation}
Z(Y_{\lambda}, T) = {{Q(\lambda, T)^{(-1)^n}}\over (1-T)(1-qT)\cdots (1-q^{n-1}T)}.
\end{equation}
By the results in \cite{W1}\cite{W2}, this one parameter family $Y_{\lambda}$ 
of toric hypersurfaces is generically ordinary for every $n$ and every prime number 
$p$.

Now, we are ready to compare the two zeta functions $Z(X_{\lambda}, T)$ and 
$Z(Y_{\lambda}, T)$. 
Let now $\lambda\in \fq$ such that $X_{\lambda}$ is smooth and $g$ is $\Delta$-regular. 
The above description shows that  
$${Z(X_{\lambda}, T) \over Z(Y_{\lambda}, T)} = ({P(\lambda, T)\over Q(\lambda, T)})^{(-1)^n}.$$
To understand this quotient of zeta functions, it suffices to understand 
the quotient ${P(\lambda, T)/Q(\lambda, T)}$.

\begin{lemma} The polynomial $Q(\lambda, T)$ divides $P(\lambda, T)$. 
\end{lemma}

{\bf Proof}. We consider the finite Galois covering $X_{\lambda} \rightarrow Y_{\lambda}$ 
with Galois group $G$, where $G=(\Z /(n+1)\Z)^{n-1}$ is an abelian group. 
For an $\ell$-adic representation $\rho: G\rightarrow {\rm GL}(V_{\rho})$, 
let $L(X_{\lambda}, \rho, T)$ 
denote the corresponding L-function of $\rho$ associated to this Galois covering. 
Then, we have the standard factorization 
$$Z(X_{\lambda}, T) = \prod_{\rho} L(X_{\lambda}, \rho, T),$$
where $\rho$ runs over all irreducible (necessarily one-dimensional) 
$\ell$-adic representations 
of $G$. If $\rho=1$ is the trivial representation, then 
$$L(X_{\lambda}, 1, T)=Z(Y_{\lambda}, T).$$
For a prime number $\ell\not=p$, 
the $\ell$-adic trace formula for $Z(X_{\lambda}, T)$ is 
$$Z(X_{\lambda}, T)= \prod_{i=0}^{2(n-1)} {\rm det}(I -T {\rm Frob}_q
|H^i(X_{\lambda}\otimes \fqbar, \Q_{\ell}))^{(-1)^{i-1}},$$
where ${\rm Frob}_q$ denotes the geometric Frobenius element over $\fq$. 
Since $X_{\lambda}$ is a smooth projective hypersurface of dimension $n-1$, one has the more 
precise form of the zeta function:
\begin{equation}
Z(X_{\lambda}, T)= {{\rm det}(I -T {\rm Frob}_q
|H^{n-1}(X_{\lambda}\otimes \fqbar, \Q_{\ell}))^{(-1)^{n}} \over 
(1-T)(1-qT)\cdots (1-q^{n-1}T)}.
\end{equation}
Similarly, the $\ell$-adic trace formula for the L-function is 
$$L(X_{\lambda}, \rho, T)= \prod_{i=0}^{2(n-1)} {\rm det}(I -T ({\rm Frob}_q\otimes 1)
|(H^i(X_{\lambda}\otimes \fqbar, \Q_{\ell}) \otimes V_{\rho})^G)^{(-1)^{i-1}}.$$
For odd $i\not=n-1$, 
$$H^i(X_{\lambda}\otimes \fqbar, \Q_{\ell})=0, \ 
(H^i(X_{\lambda}\otimes \fqbar, \Q_{\ell}) \otimes V_{\rho})^G=0.$$ 
For even $i=2k\not=n-1$ with $0\leq k\leq n-1$, 
$$H^{2k}(X_{\lambda}\otimes \fqbar, \Q_{\ell})=
\Q_{\ell}(-k), \  
(H^{2k}(X_{\lambda}\otimes \fqbar, \Q_{\ell}) \otimes V_{\rho})^G=0$$ 
for non-trivial irreducible $\rho$. 
This proves that for irreducible $\rho\not=1$, we have 
$$L(X_{\lambda}, \rho, T)= {\rm det}(I -T ({\rm Frob}_q\otimes 1)
|(H^{n-1}(X_{\lambda}\otimes \fqbar, \Q_{\ell}) \otimes V_{\rho})^G)^{(-1)^{n}}.$$
Similarly, taking $\rho=1$, one finds that 
\begin{equation}
Z(Y_{\lambda}, T)= {{\rm det}(I -T {\rm Frob}_q
|(H^{n-1}(X_{\lambda}\otimes \fqbar, \Q_{\ell}))^G)^{(-1)^{n}} \over 
(1-T)(1-qT)\cdots (1-q^{n-1}T)}.
\end{equation}
Comparing (13)-(16), we conclude that 
$$P(\lambda, T)= {\rm det}(I -T {\rm Frob}_q
|H^{n-1}(X_{\lambda}\otimes \fqbar, \Q_{\ell})),$$
$$Q(\lambda, T)= {\rm det}(I -T {\rm Frob}_q
|(H^{n-1}(X_{\lambda}\otimes \fqbar, \Q_{\ell}))^G).$$
Furthermore, the quotient 
$${P(\lambda, T) \over Q(\lambda, T)} =\prod_{\rho\not=1} 
{\rm det}(I -T ({\rm Frob}_q\otimes 1)
|(H^{n-1}(X_{\lambda}\otimes \fqbar, \Q_{\ell}) \otimes V_{\rho})^G)$$
is a polynomial with integer coefficients of degree ${n(n^n-(-1)^n)\over n+1} -n$, 
pure of weight $n-1$. The lemma is proved. 

This lemma together with Theorem 6.1 gives the following result.

\begin{theorem}There is a polynomial $R_n(\lambda, T)\in 1+T\Z [T]$ 
which is pure of weight $n-3$ and of degree ${n(n^n-(-1)^n)\over n+1} -n$,  
such that   
$${P(\lambda, T)\over Q(\lambda, T)} =R_n(\lambda, qT).$$
\end{theorem}

The polynomial $R_n(\lambda, T)$ measures how far the zeta function of 
$Y_{\lambda}$ 
differs from the zeta function of  $X_{\lambda}$. Being of integral pure weight $n-3$, 
the polynomial $R_n(\lambda, T)$  
should come from the zeta function 
of a variety (or motive $M_n(\lambda)$) of dimension $n-3$. It would be interesting 
to find this variety or motive $M_n(\lambda)$ parametrized by $\lambda$. 
In this direction, the following is known. 

If $n=2$, then $n-3<0$, $M_2(\lambda)$ is empty and we have $R_2(\lambda, T)=1$. 
If $n=3$, then $n-3=0$ and 
$$R_3(\lambda, T)=\prod_{i=1}^{18} (1-\alpha_i(\lambda)T)$$
is a polynomial of degree $18$ with $\alpha_i(\lambda)$ being roots of unity. 
In fact, Dwork \cite{Dw} proved 
that all $\alpha_i(\lambda)=\pm 1$ in this case.
Thus, $R_3(\lambda, T)$ comes from the 
the zeta function of a zero-dimensional variety $M_3(\lambda)$ 
parameterized by $\lambda$. 
What is this zero-dimensional variety $M_3(\lambda)$? 
For every $p$ and generic $\lambda$, 
the slope zeta function has the form 
$S_p(Y_{\lambda}, u, T)=1$ and 
$$S_p(X_{\lambda}, u, T)={1\over (1-T)^2(1-uT)^{20}(1-u^2T)^2}.$$
Note that $Y_{\lambda}$ is singular and not a smooth 
mirror of $X_{\lambda}$ yet. Thus, it is not surprising that the two slope zeta functions 
$S_p(X_{\lambda}, u, T)$ and $S_p(Y_{\lambda}, u, T)$ do not 
satisfy the expected slope mirror symmetry.

If $n=4$, then $n-3=1$ and 
$$R_4(\lambda, T)=\prod_{i=1}^{200}(1-\alpha_i(\lambda)T)$$
is a polynomial of degree $200$ with $\alpha_i(\lambda)=\sqrt{q}$. 
Thus, $M_4(\lambda)$ should come from  
some curve parameterized by $\lambda$. This curves has been constructed 
explicitly in a recent paper by Candelas, de la Ossa and 
Fernando-Rodriquez \cite{Ca}.   
For every $p$ and generic $\lambda$, we know that 
$S_p(Y_{\lambda}, u, T)=1$, but  
as indicated at the beginning of this 
section, we do not know if 
the slope zeta function of $X_{\lambda}$ for a generic $\lambda$ has the form 
$$S_p(X_{\lambda}, u, T)={(1-T)(1-uT)^{101}(1-u^2T)^{101}(1-u^3T) 
\over (1-T)(1-uT)(1-u^2T)(1-u^3T)}.$$ 

For general $n$ and $\lambda \in K$ for some field $K$, 
in terms of $\ell$-adic Galois representations, the pure motive 
$M_n(\lambda)$ is simply given by 
$$M_n(\lambda) = (\bigoplus_{\rho\not=1} 
(H^{n-1}(X_{\lambda}\otimes {\bar K}, \Q_{\ell}) \otimes V_{\rho})^G) \otimes 
\Q_{\ell}(-1),$$
where $\Q_{\ell}(-1)$ denotes the Tate twist. 
If $\lambda$ is in a number field $K$, this implies that the Hasse-Weil 
zeta functions of $X_{\lambda}$ and $Y_{\lambda}$ are related by 
$$\zeta(X_{\lambda}, s) = \zeta(Y_{\lambda}, s) L(M_n(\lambda), s-1).$$

\section{Slope zeta functions}

The slope zeta function satisfies a functional equation. 
This follows from the usual functional equation 
which in turn is a consequence of the Poincare duality 
for $\ell$-adic cohomology. 

\begin{proposition} 
Let $X$ be a connected 
smooth projective variety of dimension $d$ over $\fq$.  
Then the slope zeta function $S_p(X, u, T)$ 
satisfies the following functional equation 
\begin{equation}
S_p(X, u, {1\over u^dT}) = S_p(X, u, T)(-u^{d/2}T)^{e(X)},
\end{equation}
where $e(X)$ denotes the the $\ell$-adic Euler characteristic of $X$. 
\end{proposition}

{\bf Proof}. Let $P_i(T)$ denote the characteristic polynomial of 
the geometric Frobenius acting on the $i$-th $\ell$-adic cohomology of $X\otimes \fqbar$. 
Then, 
$$Z(X, T) =\prod_{i=0}^{2d}P_i(T)^{(-1)^{i+1}}.$$
Let $s_{ij}$ ($j=1,\cdots, b_i$) denote the slopes of the 
polynomial $P_i(T)$, where $b_i$ is the degree of $P_i(T)$ 
which is the $i$-th Betti number. Write
$$Q_i(T)=\prod_{j=1}^{b_i}(1-u^{s_{ij}}T).$$
Then, by the definition of the slope zeta function, we have 
$$S_p(X,u,T) = \prod_{i=0}^{2d}Q_i(T)^{(-1)^{i+1}}.$$
For each $0\leq i\leq 2d$, the slopes of $P_i(T)$ satisfies the determinant relation 
$$\sum_{j=1}^{b_i} s_{ij} = {i\over 2}b_i.$$
Using this, one computes that 
$$Q_i({1\over T}) = (-1/T)^{b_i} u^{ib_i/2}\prod_{j=1}^{b_i}(1-u^{-s_{ij}}T).$$
Replacing $T$ by $u^dT$, we get 
$$Q_i({1\over u^dT}) = ({-1\over u^dT})^{b_i} 
u^{ib_i/2}\prod_{j=1}^{b_i}(1-u^{d-s_{ij}}T).$$
The functional equation for the usual zeta function $Z(X,T)$ implies that 
$d-s_{ij}$ ($j=1,\cdots, b_i$) are exactly the slopes for $P_{2d-i}(T)$. 
Thus, 
$$Q_i({1\over u^d T}) = ({-1\over u^dT})^{b_i} u^{ib_i/2}Q_{2d-i}(T).$$
We deduce that 
$$S_p(X, u, {1\over u^dT}) = \prod_{i=0}^{2d}
\big( Q_{2d-i}(T)({-1\over u^dT})^{b_i}u^{ib_i/2}\big)^{(-1)^{i+1}}.$$
Since $b_i =b_{2d-i}$, it is clear that 
$$\sum_{i=0}^{2d}(-1)^{i}{i\over 2}b_i ={d \over 2}e(X).$$
We conclude that 
$$S_p(X,u,{1 \over u^dT})= S_p(X, u,T) (-T)^{e(X)} u^{{d\over 2}e(X)}.$$
The proposition is proved. 

From now on, we assume that $X$ is a smooth projective scheme over $W(\fq)$. 
Assume that the reduction $X\otimes \fq$ is ordinary, i.e., 
the $p$-adic Newton polygon coincides with the Hodge polygon \cite{Ma}. 
This means that the slopes of $P_i(T)$ are 
exactly $j$ ($0\leq j\leq i$) with multiplicity $h^{j, i-j}(X)$. In this case,  
one gets the explicit formula 
\begin{equation}
S_p(X\otimes \fq, u, T) = \prod_{j=0}^{d} (1-u^jT)^{e_j(X)},
\end{equation}
where
\begin{equation}
e_j(X) = (-1)^j \sum_{i=0}^d (-1)^{i-1}h^{j,i}(X).
\end{equation}
If $X$ and $Y$ form a mirror pair over the Witt ring $W(\fq)$, 
the Hodge symmetry $h^{j,i}(X)=h^{j,d-i}(Y)$ 
implies for each $j$, 
$$e_j(X)=(-1)^j\sum_{i=0}^d (-1)^{i-1}h^{j,d-i}(Y) =
(-1)^d e_j(Y).$$ 
We obtain the following result. 

\begin{proposition} Let $X$ and $Y$ be a mirror pair of $d$-dimensional smooth projective 
Calabi-Yau schemes over $W(\fq)$. Assume that both $X\otimes \fq$ and $Y\otimes \fq$ 
are ordinary. Then, we have the following symmetry for the slope zeta function:
$$S_p(X\otimes \fq, u, T) = S_p(Y\otimes \fq, u, T)^{(-1)^d}.$$
\end{proposition}

The converse of this proposition may not be always true. The slope mirror  
conjecture follows from the following slightly stronger 

\begin{conjecture}[Generically ordinary conjecture] Let $d\leq 3$. 
Suppose that $\{ X, Y\}$ form a maximally generic mirror pair of $d$-dimensional 
smooth projective Calabi-Yau schemes over $W(\fq)$. Then, both $X\otimes \fq$ 
and $Y\otimes \fq$ are generically ordinary. 
\end{conjecture}

For $d\leq 3$, this conjecture can be proved in the toric hypersurface case 
using the results in \cite{W1}\cite{W2}. 
For $d\geq 4$, we expect that the same conjecture holds if $p\equiv 1~({\rm mod}~D)$ 
for some positive integer $D$. This is again provable in the toric hypersurface case 
using the results in \cite{W1}. But we do not know if we can always take $D=1$, 
even in the toric hypersurface case if $d\geq 4$.

\end{document}